%
%

\magnification=\magstephalf
\input amstex
\documentstyle{amsppt}
\pageheight{9truein}
\pagewidth{6.5truein}
 
\TagsOnRight
\loadbold
 
\define \gln {\text{\rm GL}_n}
\define \gl {\text{\rm GL}}
\define \hofF {{\Cal H}(F)}
\define \mofF {{\frak m}(F)}
\define \MofF {{\Cal M}(F)}

\define \br {\Bbb R}
\define \be {\Bbb E}
\define \bc {\Bbb C}
\define \bz {\Bbb Z}
\define \bq {\Bbb Q}
\define \ux {\bold x}
\define \uX {\bold X}
\define \uY {\bold Y}
\define \uo {\bold 0}
\define \um {\bold m}
\define \ue {\bold e}
\define \ua {\bold a}
\define \uj {\bold j}
\define \uy {\bold y}
\define \uz {\bold z}
\define \uL {\bold L}
\define \uK {\bold K}
 
\topmatter
\title
Asymptotic Estimates for the Number of Integer Solutions to
Decomposable Form Inequalities 
\endtitle
 
\rightheadtext{Asymptotic Estimates}

\author
Jeffrey Lin Thunder
\endauthor
\address Department of Mathematics, Northern Illinois University,
DeKalb, IL 60115\endaddress
\email jthunder\@math.niu.edu\endemail
\thanks Research partially supported by NSF grant DMS-0100791\endthanks

\abstract
For homogeneous decomposable forms $F(\uX )$ in $n$ variables
with integer coefficients, we consider
the number of integer solutions $\ux\in\bz ^n$ to the inequality 
$|F(\ux )|\le m$ as $m\rightarrow\infty$. We give asymptotic estimates
which improve on those given previously by the author in [T1]. 
Here our error terms display desirable behaviour as a function
of the height whenever the degree of the form and the number of
variables are relatively prime.
\endabstract

\endtopmatter
 
\document
\baselineskip=20pt
 
\head Introduction\endhead
In this paper we consider homogeneous polynomials $F(\uX )$ in $n>1$ variables
with integer coefficients
which factor completely into a product of linear terms over $\bc$. Such
polynomials are called decomposable forms. We are concerned here with
the integer solutions to the Diophantine inequality
$$|F(\ux )|\le m\tag 1$$
Let $V(F)$ denote the $n$-dimensional volume of the set
of all real solutions $\ux\in\br ^n$ to the inequality $|F(\ux )|\le 1$,
so that by homogeneity $m^{n/d}V(F)$ is the measure of the set of 
$\ux\in\br ^n$
which satisfy (1). 
Denote the number of
integral solutions to (1) by $N_F(m)$.  

In a previous paper [T1] we answered several open questions regarding (1).
For example, $N_F(m)$ is finite
for all $m$
if and only if $F$ is of {\it finite type}:
$V(F)$ is finite, and the same is true for
$F$ restricted to any non-trivial subspace defined over $\bq$. 
Also proven in [T1] was the
following asymptotic estimate.

\proclaim{[T1, Theorem 3]} Let $F$ be decomposable form of degree $d$ in
$n$ variables  with integer
coefficients. If $F$ is of finite type, then
there are $a(F),\ c(F)\in\bq$ satisfying
$1\le a(F)\le{d\over n}-{1\over n(n-1)}$
and 
${(d-n)\over d}\le c(F)< {d\choose n}(d-n+1)$
such that
$$|N_F(m)-m^{n/d}V(F)|\ll m^{n-1\over d-a(F)}(1+\log m)^{n-2}
\hofF ^{c(F)},$$
where the implicit constant depends only on $n$ and $d$.
In particular, 
$$|N_F(m)-m^{n/d}V(F)|\ll m^{n(n-1)^2\over d(n-1)^2+1}(1+\log m)^{n-2}
\hofF^{{d\choose n}(d-n+1)}.$$
\endproclaim

The quantity $\hofF$ appearing here is defined as follows.
Write $F(\uX )=\prod _{i=1}^dL_i(\uX )$ where
the $L_i(\uX )\in \bc [\uX ]$
are linear forms in $n$ variables. 
Denote the coefficient vector of $L_i(\uX )$
by $\uL _i$ and let $\Vert\cdot \Vert$ denote the $L^2$ norm.
Then
$$\hofF =\prod _{i=1}^d\Vert\uL _i\Vert .$$

It is useful to note how the quantities $N_F(m),$ $V(F)$ and $\hofF$ vary with the form $F$. In this regard, an
important concept is the notion of {\it equivalent forms}. Two forms $F,G\in\bz [\uX ]$
are said to be equivalent if $F=G\circ T$ for some $T\in\gln (\bz )$. This is useful since
the quantities $N_F(m)$ and $V(F)$ are clearly unchanged when $F$ is replaced by an 
equivalent form. On the other hand, the height $\hofF$ is certainly not such a quantity. With
this in mind, we define
$$\MofF =\inf _{T\in\gln (\bz )}\{ {\Cal H}(F\circ T)\}.$$
One may then replace the $\hofF$ occuring in the theorem above with $\MofF$.
In a subsequent paper [T2] we showed how the main term in the estimate above, $m^{n/d}V(F)$,
is dependent on $\MofF$.

\proclaim{[T2, Theorem 2]} Let $F(\uX )\in \bz [\uX ]$ be a decomposable form
of degree $d$ in $n$ variables which doesn't vanish on $\bz ^n\setminus
\{\uo \}$. Suppose 
$V(F)$ is finite. Then
$$\MofF ^{-n/d}\ll V(F)\ll \MofF ^{-1/d}\big (1+(\log \MofF )^{n-1}\big
),$$
where the implicit constants depend only on $n$ and $d$.
\endproclaim

Note that any form $F$ of finite type satisfies the hypotheses of this theorem. 
This result points out a weakness in the asymptotic estimate above.
To wit, $N_F(m)$ is estimated by a quantity $m^{n/d}V(F)$ which decreases as
$\MofF$ increases, exactly opposite the behaviour of the error term of the
estimate.

Ideally, one would like an asymptotic estimate for the number of solutions $N_F(m)$ which
could be usefully applied uniformly for all forms $F$ of finite type, i.e., where
the error term is always dominated by the estimate $m^{n/d}V(F)$. 
Unfortunately, such can't be the case. For
example, suppose $F$ is a binary form of the kind $F(X,Y)=X^d+ \cdots$.
Then $N_F(m)\ge 2[m^{1/d}]$, where $[\ \cdot\  ]$ is the
greatest integer function. But $\MofF$ cannot
be bounded above (and whence $V(F)$ cannot be bounded below) merely because the leading coefficient
is 1. In general, simply knowing $\MofF$ is large doesn't rule out
the possibility of a great many (roughly $m^{(n-1)/d}$) integer solutions
to (1) lying in an $(n-1)$-dimensional subspace.

Our goal here is to improve the dependence on $F$ in the error term. Specifically, we aim to
derive an error term which, in so much as possible, decreases
as the ``height" of $F$ increases. 
To accomplish this, we
introduce the following more geometric ``height", one which has no
arithmetic encumberences and which is closely connected to the
volume $V(F)$. 
Define
$$\mofF =\inf \{ {\Cal H}(F\circ T)\},$$
where the infimum is over all $T\in\gln (\br )$ with $|\det (T)|=1$.
In [T1] the quantity $a(F)$ plays an important role.
Here we use a quantity
$a'(F)\ge a(F)$
which will play an analogous role. 
Like $a(F)$, the precise definition
of $a'(F)$ is somewhat complicated (we give the definition after Lemma 4
below). For the present, we simply note that it satisfies the
same inequalities as $a(F)$: $1\le a'(F)\le d/n$, and if it's less than
$d/n$, then $a'(F)\le {d\over n}-{1\over n(n-1)}$. Further, $a'(F)<d/n$
if $n$ and $d$ are relatively prime and $V(F)$ is finite.

By [T1, Proposition], $V(F)$ is finite if  $a(F)<d/n$.
It turns out that $V(F)$ is controlled
by the height $\mofF$ if $a'(F)<d/n$. Further, we can estimate $N_F(m)$ 
more precisely when $a'(F)<d/n$. 

\proclaim{Theorem 1} Let $F(\uX )\in \br [\uX ]$ be a decomposable
form of degree $d$ in $n$ variables and suppose $V(F)$ is finite.
Then $\mofF$ is an attained positive minimum and
$V(F)\ge (2/n)^n \mofF ^{-n/d}.$ 
If $F(\uX )\in\bz [\uX ]$ and $F$ doesn't vanish at any non-trivial rational
point, then $\mofF \ge n^{-d(n+1/2)/n}$.
If $a'(F)<d/n$ (in particular, if $n$ and $d$ are relatively prime), 
then $V(F)\ll\mofF ^{-n/d},$ where the implicit constant
depends only on $n$ and $d$.
\endproclaim

By Theorem 1, one doesn't expect as many
integer solutions to (1) when $\mofF$ is large in terms of $m$.
Specifically, if $a'(F)<d/n$ and $\mofF\ge m^{1/n},$ then 
$m^{n/d}V(F)\ll m^{(n-1)/d}.$  Yet it's possible to have $\gg m^{(n-1)/d}$ solutions in
a proper subspace, hence  one
can only expect a useful asymptotic estimate when $m^{n/d}V(F)$ is larger than
$m^{(n-1)/d}$; in particular, when $\mofF\ge m^{1/n}$ (if $a'(F)< d/n$). 
This explains the hypotheses for our asymptotic estimate of
$N_F(m)$ below.

\proclaim{Theorem 2} Let $F(\uX )\in \bz [\uX ]$ be a decomposable
form of finite type of degree $d$ in $n$ variables and suppose
$\mofF\le m^{1/n}$. If $a'(F)<d/n$ (in particular, if $n$ and $d$
are relatively prime) then $1\le a'(F)\le {d\over n}-{1\over n(n-1)}$
and
$$|N_F(m)-m^{n/d}V(F)|\ll
\left ({m\over \mofF ^{na'(F)/d}}\right )^{n-1\over d-a'(F)}
(1+\log m)^{n-2}.$$
The implicit constant here depends only on $n$ and $d$.
In particular,
$$|N_F(m)-m^{n/d}V(F)|\ll
\left ({m\over \mofF ^{1-{1\over d(n-1)}}}\right )^{n(n-1)^2\over d(n-1)^2+1}
(1+\log m)^{n-2}.$$
\endproclaim

We note that the main term is (almost) larger than the error term in Theorem 2
when $\mofF \le m^{1/n}$ (the ``almost" being due to the logarithmic term).
We can improve our estimate for $N_F(m)$  when $\mofF$ is close to or larger than
$m^{1/n}$ by
abandoning our goal of an asymptotic one and instead striving for a simple
upper bound.

\proclaim{Theorem 3} Let $F(\uX )\in\bz [\uX ]$ be a decomposable form
of degree $d$ in $n$ variables of finite type. If $a'(F)<d/n$ (in
particular, if $n$ and $d$ are relatively prime), then
$$N_F(m)\ll \left ({m\over\mofF }\right )^{n/d}+m^{(n-1)/d},$$
where the implicit constant depends only on $n$ and $d$.
\endproclaim

By Theorem 3, $N_F(m)\ll m^{(n-1)/d}$ if $\mofF\ge m^{1/n}$. When $a'(F)<d/n$, this improves on
[T2, Theorem 4], and (up to the implicit constants) [G, Theorem 2]
in the case $n=2$; it is the best one can say in general. 

Combining theorems 1-3 gives the following asymptotic estimate.
\proclaim{Corollary} Let $F(\uX )\in \bz [\uX ]$ be a decomposable form of degree $d$ in
$n$ variables of finite type. If $n$ and $d$ are relatively prime, then
$$|N_F(m)-m^{n/d}V(F)|\ll m^{n\over d+1/(n-1)^2}(1+\log m)^{n-2},$$
where the implicit constant depends only on $n$ and $d$.
\endproclaim

Note that we have $N_F(m)\ll m^{n/d}V(F) +m^{(n-1)/d}$ for $F(\uX )\in\bz [\uX ]$ of
finite type, provided that $a'(F)<d/n$. We would like to remove the hypothesis on $a'(F)$ here.
Unfortunately,
the volume $V(F)$ is not controlled solely by $\mofF$ in general;
see the examples in section 5 below. However, a simple modification of
the proof of [T2, Theorem 1] yields

\proclaim{Theorem 4} Let $F(\uX )\in\bz [\uX ]$ be a decomposable form of degree $d$ in
$n$ variables. 
If $V(F)$ is finite and $F$ doesn't vanish at any non-trivial rational point, then
$$V(F)\ll\mofF ^{-n/d}(1+|\log\mofF | )^{n-1},$$
where the implicit constant depends only on $n$ and $d$.
\endproclaim

In view of Theorem 1, this result would only be useful in the case $a'(F)=d/n$. 

\proclaim{Theorem 5}
Let $F(\uX )\in\bz [\uX ]$ be a decomposable form
of degree $d$ in $n$ variables of finite type.  Suppose $a'(F)=d/n$ and $\mofF \ge 1$.
If 
$$\mofF ^{-n/d}(1+\log\mofF )^{n-1}\le m^{-1/d},$$
then
$$N_F(m)\ll m^{(n-1)/d}.$$
The implicit constant here depends only on $n$ and $d$.
\endproclaim

Theorem 5 sharpens [T2, Theorem 4]
and  [G Theorem 2] (for
the $n=2$ case). In these results, the hypothesis was
$\MofF ^{1-\varepsilon }\ge m^n$ for some positive $\varepsilon$ and the
conclusion was that the solutions are contained in $\ll \varepsilon ^{1-n}$
proper subspaces. 
We'll show that $\MofF\ll\mofF ^n$ (see Lemma 7), thus Theorem 5 represents a
true improvement on these results. A reasonable conjecture, in view of
our results here, is that $N_F(m)\ll m^{(n-1)/d}$ whenever $V(F)\ll m^{-1/d}$ and
$F$ is of finite type.

It is possible to explicitly determine bounds for the implicit
constants in the above results. Frankly, they wouldn't be very ``good", as
our proofs ultimately rely on quantitative versions of the subspace theorem.
We've attempted to keep some track of constants depending on $n$ and $d$
in the
following two sections and, to the extent where relatively painless, in the proofs of
our theorems. In general, very little effort has been expended trying to
get good bounds for these constants. For the remainder of this paper,
all implicit constants depend only on $n$ and $d$.

\head 1. Preparatory Lemmas 
\endhead

Throughout this section, $F(\uX )\in\br [\uX ]$ is assumed to be a
decomposable form of degree $d$ in $n$ variables. Also, all vectors
are assumed to be row vectors.

\proclaim{Lemma 1} 
Let $a>0$ and $T\in\gln (\br )$.
Then 
$${\Cal H}(aF)=a \hofF,\qquad {\frak m}(aF)=a \mofF,\qquad \text{and}\qquad {\frak m}(F\circ T)
= |\det (T)|^{d/n}\mofF .$$
Further, if $V(F)$ is finite, then
$$V(aF)=a^{-n/d}V(F)\qquad\text{and}\qquad V(F\circ T)=|\det (T)|^{-1}V(F).$$
\endproclaim

\demo{Proof} The first two equations are clear from the definitions. As for the third,
write
$T=D S,$
where $D$ is the diagonal matrix with entries $a=|\det (T)|^{1/n}$ and $S\in\gln (\br )$ with $|\det (S)|=1.$
Then $F\circ D=a^dF$ and ${\frak m}(F\circ T)={\frak m}(a^dF\circ S)={\frak m}(a^dF)=a^d\mofF$. For the
last equation, let ${\Cal S}$ be the set of all $\ux\in\br ^n$ with $|F(\ux )|\le 1,$ so that
$V(F)$ is the volume of ${\Cal S}$. Then $V(F\circ T)$ is the volume of $T^{-1}({\Cal S})$, which is
$|\det (T)|^{-1}V(F).$ Finally, the fourth equation can be viewed as a special case of the last
(write $T$ as we did above), or as a simple consequence of the homogeneity of $F$.
\enddemo

\proclaim{Lemma 2} Suppose $V(F)$ is finite. Then $\mofF$ is an attained 
positive minimum and $V(F) \ge (2/n)^n\mofF ^{-n/d}$.
\endproclaim

\demo{Proof} 
For a $T\in\gln (\br )$ with $|\det T|=1$, let $P(T)$ be the parallelepiped defined by
$$P(T)=\{a_1\ux _1+\cdots +a_n\ux _n\: |a_j|\le 1\ \text{for all $1\le j\le n$}\},$$
where $\ux _1^{tr},\ldots ,\ux _n ^{tr}$ are the columns of $T$.
Note that the volume of $P(T)$ is $2^n$ and that
$$\prod _{i=1}^d\max _{1\le j\le n}\{|L_i(\ux _j)|\}\le {\Cal H}(F\circ T) .$$
In particular, $|F(\ux )|\le n^d{\Cal H}(F\circ T)$ for all $\ux$ in $P(T)$.

Let ${\Cal C}$ be the set of all $T\in\gln (\br)$ with $|\det (T)|=1$ and
${\Cal H}(F\circ T)\le 2\mofF$. 
Suppose ${\Cal C}$ is unbounded (viewed as a subset of $\br ^{n^2}$ in the usual way). Then for
some $1\le j_0\le n$ there is an infinite sequence $T_1,\ldots \in {\Cal C}$ where, letting
$\ux _{i,1}^{tr},\ldots ,\ux _{i,n}^{tr}$ denote the columns of $T_i$,  
we have $\Vert \ux _{i+1,j_0}\Vert
\ge 2\Vert\ux _{i,j_0}\Vert$ for all $i\ge 1$. But this implies the existence of an infinite
sequence of parallelepipeds $P(T_1),\ldots$, all of which are contained in the set
$$\{\ux\in\br ^n\: |F(\ux )|\le n^d2\mofF\},$$
and also satisfying
$$Vol \left (P(T_{i+1})\setminus\cup _{l\le i}P(T_l)\right )\ge 2^{-1}Vol\left (P(T_{i+1})\right )
=2^{n-1}.$$ This contradicts the hypothesis that $V(F)$ is finite, thus ${\Cal C}$ is
bounded. 

The map $T\mapsto {\Cal H}(F\circ T)$ is clearly continuous. Since ${\Cal C}$ is bounded (and
certainly closed), $\mofF$ is an attained minimum. Let $\mofF ={\Cal H}(F\circ T)$ for some
$T\in {\Cal C}$. Then 
$|F(\ux )|\le n^d\mofF$ for all $\ux$ in $P(T)$,
which implies that $V\big (n^{-d}\mofF ^{-1}F\big )\ge 2^n$ and 
$V(F) \ge (2/n)^n\mofF ^{-n/d}$ by Lemma 1.
\enddemo

\proclaim{Lemma 3} Suppose $V(F)$ is finite. Write
$F(\uX )=\prod _{i=1}^dL_i(\uX)$ where the $L_i(\uX )$ are real linear
forms for $i\le r$ and complex for $i=r+1,\ldots ,d=r+2s$, with
$\uL _{i+s}=\overline{\uL _i}$ for $i=r+1,\ldots ,r+s$. Suppose
$a_1,\ldots ,a_d$ are positive real numbers whose product is 1, so that
$F(\uX )=\prod _{i=1}^da_iL_i(\uX )$. Then

$$\sum _{1\le i_1,\ldots ,i_n\le d}|\det (a_{i_1}\uL_{i_1}^{tr}\cdots 
a_{i_n}\uL _{i_n}^{tr})|^2\ge {n!\over n^n}d^n\mofF ^{2n/d}.$$
\endproclaim

\demo{Proof}
Let $\sigma$ be the permutation of $\{1,\ldots ,d\}$ induced by complex
conjugation, i.e., 
$$\sigma (i)=\cases i&\text{if $i\le r$,}\\
i+s&\text{if $r<i\le r+s$},\\
i-s&\text{if $r+s<i\le d$}.\endcases$$
Let $b_i$ be the geometric mean of $a_i$ and $a_{\sigma (i)}$. Then
the product $\prod _{i=1}^db_i=\prod _{i=1}^da_i$ and
for any $n$-tuple $(i_1,\ldots ,i_n)$ we have
$(a_{i_1}\cdots a_{i_n})^2+(a_{\sigma (i_1)}\cdots a_{\sigma (i_n)})^2\ge
2(b_{i_1}\cdots b_{i_n})^2.$ Since 
$$|\det (\uL_{i_1}^{tr}\cdots 
\uL _{i_n}^{tr})|^2=
|\det (\uL_{\sigma (i_1)}^{tr}\cdots 
\uL _{\sigma (i_n)}^{tr})|^2,$$
we see that
$$\aligned 2\sum _{1\le i_1,\ldots ,i_n\le d}|\det (a_{i_1}\uL_{i_1}^{tr}\cdots 
a_{i_n}\uL _{i_n}^{tr})|^2&=
\sum _{1\le i_1,\ldots ,i_n\le d}(a_{i_1}\cdots a_{i_n})^2|\det (\uL_{i_1}^{tr}
\cdots \uL _{i_n}^{tr})|^2\\
&\qquad +
\sum _{1\le i_1,\ldots ,i_n\le d}(a_{\sigma (i_1)}\cdots a_{\sigma (i_n)})^2
|\det (\uL_{i_1}^{tr}\cdots \uL _{i_n}^{tr})|^2\\
&\ge 2\sum _{1\le i_1,\ldots ,i_n\le d}|\det (b_{i_1}\uL_{i_1}^{tr}\cdots 
b_{i_n}\uL _{i_n}^{tr})|^2.\endaligned$$
The upshot is that we may replace $a_i$ and $a_{\sigma (i)}$ with $b_i$, i.e.,
we may assume $a_i=a_{\sigma (i)}$.
But if this is the case, then 
it suffices to prove the lemma under the assumption that
$a_i=1$ for
all $i$. 

Let $\be ^d\subset\br ^{r}\oplus 
\bc ^{2s}$ 
be the set of
$\ux =(x_1,\ldots ,x_d)$ where $x_{i+s}=\overline{x_i}$ for $r+1\le i\le r+s$.
Then $\be ^d$ is  $d$-dimensional Euclidean space via the usual
hermitian inner product on $\bc ^d$. 
If $M$ is the $d\times n$ matrix with rows $\uL _1,\ldots ,\uL _d$, then
the columns of $M$ are in $\be ^d$. Moreover, the rank of $M$ must
be $n$ since $V(F)$ is finite (see [T2]).
We apply Gram-Schmidt to the matrix $M$;
there is an upper triangular $T\in\gln (\br)$ such that $MT$ is a matrix
with orthonormal columns (in $\be ^d$). Denote the rows of $MT$ by
$\uL _1',\ldots ,\uL _d'$ and the columns by $\um _1^{tr},\ldots ,\um _n^{tr}$.

Using the inequality between the arithmetic and geometric means and also
Lemma 1, we have
$$\aligned
{n\over d}={1\over d}\sum _{j=1}^n\Vert \um _j\Vert ^2&={1\over d}\sum 
_{i=1}^d\Vert\uL _i'\Vert ^2\\
&={1\over d}\sum _{i=1}^d\Vert \uL _iT\Vert ^2\\
&\ge \left (\prod _{i=1}^d\Vert \uL _iT\Vert ^2\right )^{1/d}\\
&={\Cal H}(F\circ T)^{2/d}\\
&\ge {\frak m}(F\circ T)^{2/d}\\
&=\mofF ^{2/d}|\det T|^{2/n}.\endaligned$$
Hence
$|\det T|^2\le \left ({n\over d}\right )^n\mofF ^{-2n/d}.$
On the other hand, since the $\um _j$s are orthonormal,
$$\aligned 1=\Vert\um _1\wedge \cdots\wedge \um _n\Vert &=
\sum _{1\le i_1<\cdots <i_n\le d}|\det \big ((\uL _{i_1}')^{tr}
\cdots
(\uL _{i_n}')^{tr}\big )|^2\\
&=
\sum _{1\le i_1<\cdots <i_n\le d}|\det (\uL _{i_1}^{tr}\cdots
\uL _{i_n}^{tr})|^2|\det T|^2\\
&\le
\sum _{1\le i_1<\cdots <i_n\le d}|\det (\uL _{i_1}^{tr}\cdots
\uL _{i_n}^{tr})|^2\left ({n\over d}\right )^n\mofF ^{-2n/d}.\endaligned$$
This inequality suffices to prove the lemma.
\enddemo

Let
$$c_1= {(d/n)^n\over {d\choose n}}.$$

\proclaim{Lemma 4} Suppose $V(F)$ is finite and $\hofF =\mofF$.
Let $A<1$ and $1\le j<n$ and suppose there is an
$S\subset \{1,\ldots ,d\}$ with cardinality
$|S|=[jd/n]+1$ such that
$${\Vert \uL _{i_1}\wedge \cdots \wedge \uL _{i_{j+1}}\Vert\over
\Vert \uL _{i_1}\Vert\cdots \Vert\uL _{i_{j+1}}\Vert }\le A$$
for all $i_1,\ldots ,i_{j+1}\in S$. Then there is a factorization
$F(\uX )=\prod _{i=1}^da_iL_i(\uX )$ as in Lemma 3 with
$$\sum _{1\le i_1,\ldots ,i_n\le d}|\det (a_{i_1}\uL_{i_1}^{tr}\cdots 
a_{i_n}\uL _{i_n}^{tr})|^2\le n!{d\choose n}A^{n([jd/n]+1)-jd\over (n-j)d}
\mofF ^{2n/d}.$$
In particular, 
$$A\ge c_1^{(n-j)d\over n([jd/n]+1)-jd}.$$ 
\endproclaim

\demo{Proof}
Without loss of generality we may assume $\mofF =1$ and $\Vert\uL _i\Vert
=1$ for all $i=1,\ldots ,d$.

Let $a=A^{[jd/n]+1-d\over (n-j)d}$, which is greater than 1 since $A<1$.
Let $b=a^{-[jd/n]-1\over d-[jd/n]-1}$, so that $b<1<a$ and
$a^{[jd/n]+1}b^{d-[jd/n]-1}=1$. We note that $a^jb^{n-j}=a^nA=A^{
n([jd/n]+1)-jd\over (n-j)d}.$

For our factorization, let $a_i=a$ if $i\in S$ and $a_i=b$ otherwise.
Given an $n$-tuple $(i_1,\ldots ,i_n)$ with $l$ of the indices in $S$,
we have
$$|\det (a_{i_1}\uL _{i_1}^{tr}\cdots a_{i_n}\uL _{i_n}^{tr})|\le a^lb^{n-l}
\le a^jb^{n-j}$$ 
if $l<j+1$, and
$$|\det (a_{i_1}\uL _{i_1}^{tr}\cdots a_{i_n}\uL _{i_n}^{tr})|\le a^lb^{n-l}
A\le a^nA$$
if $l\ge j+1$. Lemma 4 follows from this and Lemma 3.
\enddemo

We can now state with more clarity exactly what the quantity $a'(F)$ is.
Suppose $\mofF =\hofF$. For $1\le j\le n-1$, let $s_j(F)$ be the 
cardinality of the largest subset $S\subset \{1,\ldots ,d\}$ where
$${\Vert \uL _{i_1}\wedge \cdots \wedge \uL _{i_{j+1}}\Vert\over
\Vert \uL _{i_1}\Vert\cdots \Vert\uL _{i_{j+1}}\Vert }
<c_1^{(n-j)d\over n([jd/n]+1)-jd}$$ 
for all $i_1,\ldots ,i_{j+1}\in S$. Note that we don't demand the
$i_j$'s be distinct, so any set with just one
element will vacuously satisfy this criterion. By Lemma 4, $s_j(F)\le [jd/n]$.
Let $s(F)$ be the maximum of $s_j(F)/j$ over all $1\le j\le n-1$.
Then $1\le s(F)\le d/n$. Moreover, $s(F)\le {d\over n}-{1\over n(n-1)}$
if $s(F)<d/n$.

For an arbitrary $F$, we define
$$a'(F):=\max \{s(F\circ T)\},$$
where the maximum is over all $T\in\gln (\br )$ with $|\det (T)|=1$
and $\mofF ={\Cal H}(F\circ T)$. Then $1\le a'(F)\le d/n$, and
$a'(F)\le {d\over n}-{1\over n(n-1)}$ if $a'(F)<d/n$. Moreover,
$a'(F)<d/n$ if $n\nmid jd$ for all $j=1,\ldots ,d$, i.e., if $n$ and
$d$ are relatively prime. 

\proclaim{Lemma 5a} Let $M,N\ge 1$ and fix $\uK _1,\ldots ,\uK _N\in\bc ^M$. Then
for any $\uL _1,\ldots ,\uL _{N+1}\in\bc ^M$ with $\Vert\uL _i\Vert =1$
for all $i$, we have
$$\sum _{j=1}^{N+1}\Vert\uK _1\wedge \cdots \wedge \uK _n\wedge\uL _j\Vert
^2\ge \Vert\uK_1\wedge\cdots \wedge\uK _N\Vert ^2\cdot\Vert\uL _1\wedge
\cdots\wedge\uL _{N+1}\Vert ^2.$$
\endproclaim

\demo{Proof} This is trivial if the $\uL_i$'s are linearly dependent,
so assume otherwise. After possibly applying a unitary transformation,
we may assume that the span of $\uL _1,\ldots \uL _j$ is equal to the
span of the first $j$ canonical basis vectors $\ue_1,\ldots ,\ue _j$
for $j=1,\ldots ,N+1$. Let $a_j=\uL _j\cdot\ue _j$ for $j=1,\ldots ,N+1$
and write $\uK _i=\sum _{l=1}^Mk_{i,l}\ue _l$ for $i=1,\ldots ,N$.

By [S, Chap. I, Lemma 5A], we see that
$$\Vert\uK _1\wedge\cdots\wedge\uK _N\wedge\uL _j\Vert ^2\ge\sum _{\sigma}
|\det \Sb 1\le i\le N\\ l\in\sigma\endSb (k_{i,l}) |^2\cdot |a_j|^2$$
for all $j=1,\ldots ,N+1$, where the sum is over all $N$-tuples
$\sigma = (i_1,\ldots ,i_N)$ with $i_1<\cdots <i_N$ and $j\not\in\sigma$.
On the other hand
$$\Vert\uK_1\wedge\cdots \wedge\uK _N\Vert ^2\cdot\Vert\uL _1\wedge
\cdots\wedge\uL _{N+1}\Vert ^2=\sum _{\sigma}|\det \Sb 1\le i\le N\\
l\in\sigma\endSb (k_{i,l}) |^2|a_1|^2\cdots |a_{N+1}|^2,$$
where the sum is over all $N$-tuples
$\sigma = (i_1,\ldots ,i_N)$ with $i_1<\cdots <i_N$. But for any
such $\sigma$, there is a $j\in \{1,\ldots ,N+1\}$ with $j\not\in\sigma$.
Further, since $\Vert\uL _j\Vert =1$ for all $j$,
$|a_j|^2\ge |a_1|^2\cdots |a_{N+1}|^2$ for all $j$. This proves the lemma.
\enddemo

\proclaim{Lemma 5b} Suppose $V(F)$ is finite and $\mofF =\hofF$. Let  $B>1$,
$1\le j<n$, and fix
linearly independent
$\uL _{i_1},\ldots ,\uL _{i_j}$. Suppose $S\subset \{ 1,\ldots ,d\}$ with
cardinality $|S|=[jd/n]+1$ such that
$${\Vert \uL _{i_1}\wedge\cdots\wedge\uL _{i_j}\wedge \uL _l\Vert\over
\Vert \uL _{i_1}\wedge\cdots\wedge\uL _{i_j}\Vert\cdot \Vert\uL _l\Vert}
\le B$$
for all $l\in S$. Then
$$\sqrt{j+1}B\ge c_1^{(n-j)d\over n([jd/n]+1)-jd}.$$ 
\endproclaim

\demo{Proof} Without loss of generality, we may assume $\mofF =1$ and
$\Vert\uL _i\Vert =1$ for all $i$. Let $l_1,\ldots ,l_{j+1}\in S$. Then
by Lemma 5a and the hypotheses,
$$\Vert\uL _{l_1}\wedge\cdots\wedge\uL _{l_{j+1}}\Vert ^2\le
\sum _{k=1}^{j+1}{\Vert \uL _{i_1}\wedge\cdots\wedge\uL _{i_j}\wedge
\uL _{l_k}\Vert ^2\over\Vert\uL _{i_1}\wedge\cdots\wedge\uL _{i_j}\Vert ^2}
\le (j+1)B^2.$$
Lemma 5b follows from this and Lemma 4.
\enddemo

We now come to our fundamental inequality. This result will be
used as an alternative to [T1 Lemma 5]. (It can also be used in
place of [T1 Lemma 6] in the case when $a'(F)=d/n$.) 

\proclaim{Lemma 6} Suppose $V(F)$ is finite and $\mofF =\hofF$. 
Then for any $\ux\in\br ^n$, 
there are $n$ linearly independent linear factors $L_{i_1}(\uX ),
\ldots ,L_{i_n}(\uX )$ of $F(\uX )$ satisfying
$$\left ({\prod _{j=1}^n|L_{i_j}(\ux )|\over|\det (\uL _{i_1}^{tr}\cdots
\uL _{i_n}^{tr})|}\right )^{a'(F)}\le c_2 {|F(\ux )|\over
\|\ux \|^{d-na'(F)}\mofF},$$
where
$$c_2= n^{n(d-na'(F))/2}\left ((n!)^{1/2}\prod _{j=1}^{n-1}c_1^{(j-n)d\over 
n([jd/n]+1)-jd}\right )^{d-(n-1)a'(F)}.$$

If $V(F)$ is finite and $T\in\gln (\br )$ satisfies $|\det T|=1$
and $\mofF = {\Cal H}(F\circ T)$, 
then for any $\ux\in\br ^n$, 
there are $n$ linearly independent linear factors $L_{i_1}(\uX ),
\ldots ,L_{i_n}(\uX )$ of $F(\uX )$ satisfying
$$\left ({\prod _{j=1}^n|L_{i_j}(\ux )|\over|\det (\uL _{i_1}^{tr}\cdots
\uL _{i_n}^{tr})|}\right )^{a'(F)}\le c_2 {|F(\ux )|\over
\|T^{-1}(\ux ) \|^{d-na'(F)}\mofF}.$$
\endproclaim

\demo{Proof} We first note how the second part follows directly from
the first. Given such a $T$ and $\ux$, 
$F\circ T \big (T^{-1}(\ux )\big )=F(\ux )$, and 
similarly for each linear factor $L_i(\ux )$. Apply the first part
of the lemma to $F\circ T$ and $T^{-1}(\ux )$ to obtain the second part.

As remarked in the proof of Lemma 3, there are $n$ linearly
independent factors of $F$ if $V(F)$ is finite. So if $F(\ux )=0$, the
lemma is trivially true. Suppose now that $F(\ux )\neq 0$. 
By homogeneity and Lemma 1, we may assume without loss of generality
that $\hofF =1$, and further that $\Vert\uL _i\Vert =1$ for all $i=1,
\ldots ,d$. For 
notational convenience, set
$$c_{1,j}={c_1^{(n-j)d\over n([jd/n]+1)-jd}\over \sqrt{j+1}}$$
for $j=1,\ldots ,n-1$.

Let
$|L_{i_1}(\ux )|=\min_{1\le i\le d}\{|L_i(\ux )|\}$
and let $S_1\subset \{1,\ldots ,d\}$ be the
subset of indices $l$ such that
$$\Vert \uL _{i_1}\wedge\uL _l\Vert <c_{1,1}.$$
Continue recursively in the following manner: for $j>1$ let
$|L_{i_j}(\ux )|$ 
be the minimum of
$|L_i(\ux )|$ over all $i$ not in $S_{j-1}$
and let $S_j\subset
\{1,\ldots ,d\}$ be the subset of indices $l$ with
$${\Vert \uL _{i_1}\wedge\cdots\wedge\uL _{i_j}\wedge \uL _l\Vert\over
\Vert \uL _{i_1}\wedge\cdots\wedge\uL _{i_j}\Vert}
<c_{1,j}.$$
By definition,
$|S_j|\le ja'(F)\le [jd/n]<d$ for each $j=1,\ldots ,n-1,$ allowing us to
continue up to a choice for $L_{i_n}(\ux )$. Note that $|S_n|=d$.

By construction, we have $|L_{i_1}(\ux )|\le \cdots \le |L_{i_n}(\ux )|$
and
$$|F(\ux)|=\prod _{i=1}^d|L_i(\ux )|\ge \prod _{j=1}^n|L_{i_j}(\ux )|
^{a_j},$$
where $a_1=|S_1|$ and $a_j=|S_j|-|S_{j-1}|$ for $j>1$. Note in particular
that $a_1+\cdots +a_j=|S_j|\le ja'(F)$ for all $j<n$. Letting $s=
(n-1)a'(F)-|S_{n-1}|,$ we have
$$|F(\ux)|\ge |L_{i_1}(\ux )|^{a_1}\cdots |L_{i_{n-2}}(\ux )|^{a_{n-2}}
\cdot |L_{i_{n-1}}(\ux )|^{a_{n-1}+s}\cdot |L_{i_n}(\ux )|^{a_n-s}.$$
Since $a_1+\cdots +a_{n-2}+a_{n-1}+s=(n-1)a'(F)$ and $a_n-s=d-(n-1)a'(F)$,
[T1, Lemma 1] implies that
$$|F(\ux)|\ge \prod _{j=1}^n|L_{i_j}(\ux )|^{a'(F)}\cdot |L_{i_n}(\ux )|
^{d-na'(F)}.$$
By [T1, Lemma 4],
$$|L_{i_n}(\ux )|\ge n^{-n/2} \Vert\ux\Vert |\det (\uL _{i_1}^{tr}\cdots
\uL _{i_n}^{tr})|,$$
so that
$$\left ({\prod _{j=1}^n|L_{i_j}(\ux )|\over|\det (\uL _{i_1}^{tr}\cdots
\uL _{i_n}^{tr})|}\right )^{a'(F)}\le {n^{n(d-na'(F))/2}|F(\ux )|\over
\|\ux \|^{d-na'(F)}|\det(\uL _{i_1}^{tr}\cdots \uL _{i_n}^{tr})|
^{d-(n-1)a'(F)}}.$$
Finally, we have
$$\aligned \Vert\uL _{i_1}\wedge\cdots\wedge\uL _{i_n}\Vert &\ge
c_{1,n-1}\Vert\uL _{i_1}\wedge\cdot\wedge\uL _{i_{n-1}}\Vert\\
&\ge c_{1,n-1}c_{1,n-2}\Vert\uL _{1_1}\wedge\cdots\wedge\uL_{i_{n-2}}\Vert\\
&\vdots\\
&\ge c_{1,n-1}\cdots c_{1,1}.\endaligned$$
The lemma follows.
\enddemo

\proclaim{Lemma 7} Suppose $F(\uX )\in\bz [\uX ]$ and that $F$ doesn't vanish
on $\bz ^n\setminus\{\uo\} .$ Let ${\Cal H}(F\circ T)=\mofF$ with
$T\in\gln (\br )$, $|\det (T)|=1.$ Then there is an $S\in\gln (\bz )$ with
${\Cal H}(F\circ S)\le n^{d(n+1/2)} \mofF ^n$ such that
$$\mofF ^{-1/d}n^{-3/2}(n!)^{-2}\Vert\uy\Vert\le\Vert T^{-1}S (\uy )\Vert\le n^{n+1/2} 
\mofF ^{(n-1)/d}\Vert\uy\Vert$$
for all $\uy\in\br ^n$.
For such an $F$, $\MofF\le n^{d(n+1/2)}\mofF ^n$, and in particular
$\mofF \ge n^{-(n+1/2)d/n}$.
\endproclaim

\demo{Proof} Let $\ux _1^{tr},\ldots ,\ux _n^{tr}$ denote the columns of $T$
and let $P(T)$ be the parallelepiped defined in the proof of Lemma 2 above.
Let $\lambda _1\le \cdots
\le \lambda _n$ be the successive minima of $P(T)$ with respect to the
integer lattice $\bz ^n$. Since the volume of $P(T)$ is $2^n$, 
Minkowski's theorem says that
$$(n!)^{-1}\le \lambda _1\cdots \lambda _n\le 1.\tag 2$$

Choose a basis $\uz _1,\ldots ,\uz _n$ for $\bz ^n$ satisfying
$\uz _i\in i\lambda _iP$ for $1\le i\le n$ and let 
$S$ be the matrix with columns $\uz _1^{tr},\ldots ,\uz _n ^{tr}.$
Write
$$T^{-1}S=\pmatrix a_{1,1}&\hdots &a_{1,n}\\
          \vdots &\ddots &\vdots \\
          a_{n,1}&\hdots &a_{n,n}\endpmatrix ,$$
so that $\sum _{i=1}^na_{i,j}\ux _i^{tr}=\uz _j^{tr}$ for all $1\le j\le n$.
In particular, since $\uz _j\in j\lambda _jP$, we have
$$|a_{i,j}|\le j\lambda _j\qquad 1\le i,j\le n .\tag 3$$
Similarly, writing $S^{-1}T=(b_{i,j})$  and using Cramer's rule, we see that
$$|b_{i,j}|\le (n-1)!\prod _{l\neq i}l\lambda _l
< n!\prod _{l\neq i}l\lambda _l,\qquad 1\le i,j\le n.\tag 3'$$

As before, write $F(\uX )=\prod _{i=1}^dL_i(\uX )$. By (3),
$$\Vert \uL _iS\Vert ^2 =\Vert \uL _iTT^{-1}S\Vert ^2\le n (n\lambda _n)^2
\Vert\uL _iT\Vert ^2,$$
which implies that 
$${\Cal H}(F\circ S)\le n^{3d/2}\lambda _n^d\mofF .\tag 4$$
Also by (3), for any $\uy\in\br ^n$,
$$\Vert T^{-1}S (\uy )\Vert ^2\le n(n\lambda _n)^2\Vert\uy\Vert ^2.\tag 5$$
Using (3') in a similar manner yields
$$\Vert T^{-1}S (\uy )\Vert ^2n (n!)^2\left (2\lambda _2\cdots n\lambda _n\right )^2\ge
\Vert\uy\Vert ^2.\tag 5'$$

As seen in the proof of Lemma 2,
$|F(\uy)|\le n^d\mofF$ for all $\uy\in P$, and by homogeneity
$|F(\uz _1)|\le (n\lambda _1)^d\mofF$. But since $F(\uX )\in \bz [\uX ]$
and since $F$ doesn't vanish on $\bz ^n\setminus\{\uo\}$, we conclude
that $|F(\uz _1)|\ge 1$. Hence $\lambda _1\ge n^{-1}\mofF ^{-1/d}$. By (2), this
implies that
$\lambda _2\cdots\lambda _n\le n\mofF ^{1/d}$ and $\lambda _n\le n^{n-1}\mofF ^{(n-1)/d}$.
Lemma 7 follows from these estimates, (4), (5), (5') and [T1, Lemma 2] 
(where it is shown that $\MofF \ge 1$ if $F(\uX )\in\bz [\uX ]$).

\enddemo

We will also need the following result from [T1]:

\proclaim{Lemma 8 [T1, Lemma 7]} Let $K _1(\uX ),\ldots ,K _n(\uX )\in\bc [\uX ]$ 
be $n$ linearly independent linear forms
in $n$ variables. Denote the corresponding coefficient vectors
by $\uK _1,\ldots ,\uK _n$. Let $A,B,C>0$ with $C>B$ and let $D>1$. 
Consider the set of $\ux\in\br ^n$ satisfying
$${\prod _{i=1}^n|K_i(\ux )|\over |\det(\uK _1^{tr}\cdots \uK _n^{tr}
)|}\le A$$
and also $B\le\|\ux\|\le C.$
If $BC^{n-1}\ge D^{n-1}n!n^{n/2}A$,
then this set
lies in the union of less than 
$$n^3\left (\log _D\big (BC^{n-1}/n!n^{n/2}A\big )\right )^{n-2} $$
convex sets of the form
$$\gather
\{\uy\in\br ^n\: |K_i'(\uy )|\le a_i\ \text{for $i=1,\ldots ,n$}\}, \\
|\det \big ((\uK _1')^{tr}\cdots (\uK _n')^{tr}\big )|=1,\tag 6 \\
\|\uK '_i\| =1\qquad i=1,\ldots ,n,\endgather$$
with 
$$\prod _{i=1}^na_i<D^nn!n^{n/2}{CA\over B}.$$
If $BC^{n-1}<D^{n-1}n!n^{n/2}A,$ then this set lies in the
union of no more than $n!$ convex sets of this form.
\endproclaim

\head 2. Intermediate Results\endhead

\proclaim{Proposition 1} Let $F(\uX )\in\br [\uX ]$ be a decomposable form of degree $d$ in
$n$ variables with $V(F)$ finite. Suppose that $\hofF =\mofF$.
Let $1\le B<C$ and $D>1$ and let $\Lambda$ be a lattice of rank $n$.
Then the $\ux\in\br ^n$ with 
$$\left ({m\over \mofF}\right ) ^{1/d}B\le\|\ux\|\le 
\left ({m\over \mofF }\right )^{1/d}C$$ 
satisfying {\rm (1)} lie in no more than
$${d\choose n}\max \{n!,\ n^3\big (\log _D(B^{d-(n-1)a'(F)\over a'(F)}C^{n-1})\big )^{n-2}\}$$ 
convex sets of the form {\rm (6)} with 
$$\prod _{i=1}^na_i< c_3n^{n/2}n!\left (
{m\over \mofF }\right )^{n/d}
D^n{C\over B^{d-(n-1)a'(F)\over a'(F)}},$$
where $c_3=\max \{1,c_2^{1/a'(F)}\}$. 
Further, such a set has volume no greater than
$$c_32^nn^{n/2}(n!)^2\left (
{m\over \mofF }\right )^{n/d}
D^n{C\over B^{d-(n-1)a'(F)\over a'(F)}}$$
and either all lattice points in such a set lie in a 
sublattice of smaller rank, or the convex set contains no more than
$${c_33^n2^{n(n-1)/2}n^{n/2}(n!)^2\over\det\Lambda}\left (
{m\over \mofF }\right )^{n/d}
D^n{C\over B^{d-(n-1)a'(F)\over a'(F)}}$$
lattice points.
\endproclaim

\demo{Proof} First assume $m=\mofF =1$. By Lemma 6, if $|F(\ux )|\le 1$  and $\|\ux\|\ge 1$,
then there are $n$ linearly independent
linear factors $L_{i_1}(\uX ),\ldots ,L_{i_n}(\uX )$ of $F(\uX )$ such that
$${\prod _{j=1}^n|L_{i_j}(\ux )|\over |\det (\uL _{i_1}^{tr}\cdots \uL _{i_n}^{tr})|}
\le c_3\|\ux\| ^{-d+na'(F)\over a'(F)}.$$
In particular, if $\|\ux\|\ge B,$ then
$${\prod _{j=1}^n|L_{i_j}(\ux )|\over |\det (\uL _{i_1}^{tr}\cdots \uL _{i_n}^{tr})|}
\le c_3B^{-d+na'(F)\over a'(F)}.\tag 7$$

We now invoke Lemma 8, using $A=c_3B^{-d+na'(F)\over a'(F)}$.
Accordingly, the $\ux\in\br ^n$ with $B\le\|\ux\|\le C$ and
satisfying (7) lie in the union of no more than
$$\multline \max \{n!,\ n^3\left (\log _D\big (BC^{n-1}/n!n^{n/2}A\big )\right )^{n-2}\}\\
=\max \{ n!,\ n^3\big (\log _D(B^{d-(n-1)a'(F)\over a'(F)}C^{n-1}/n!n^{n/2}c_3)\big )^{n-2}\}\\
\le
\max \{ n!,\ n^3\big (\log _D(B^{d-(n-1)a'(F)\over a'(F)}C^{n-1})\big )^{n-2}\}\\
\endmultline$$ 
convex sets of the form (6) with 
$$\prod _{i=1}^na_i<n!n^{n/2}D^n{CA\over B}= c_3n!n^{n/2} D^n{C\over B^{d-(n-1)a'(F)\over a'(F)}}.$$
By [T1, Lemma 9], the volume of such a convex set is  no greater than $2^nn!\prod _{i=1}^na_i$ and the
number of lattice points in such a set is no more than ${3^n2^{n(n-1)/2}n!\over\det (\Lambda )}\prod _{i=1}^na_i$
if there are $n$ linearly independent lattice points in the set.

Since there are at most ${d\choose n}$ possible $n$-tuples $L_{i_1},\ldots ,L_{i_n}$
to consider above, this proves Proposition 1 when $m=\mofF =1$. For the general case, let
$G=\mofF ^{-1}F$. Then ${\frak m}(G)=1$ by Lemma 1 and $|F(\ux )|\le m$ if and
only if $|G(\ux )|\le {m\over\mofF }$. But this it true if and only if
$|G(\uy )|\le 1$, where $\ux =\left ({m\over\mofF}\right )^{1/d}\uy$. In this way, we see
that the general case follows from the case $m=\mofF=1$ via a dilation by
$\left ({m\over \mofF}\right )^{1/d}$.

\enddemo

\proclaim{Proposition 2}  Let $F(\uX )\in\br [\uX ]$ be a decomposable
form of degree $d$ in $n$ variables with $V(F)$ finite. Suppose $a'(F)<d/n$
and $\hofF =\mofF$. Let $B_0\ge 1$ and $D>1$, and let $\Lambda$ be
a lattice of rank $n$. Then the volume of all $\ux$ 
satisfying {\rm (1)} and
with $\|\ux\|\ge B_0(m/\mofF )^{1/d}$ is less than
$$c_4\left ({m\over\mofF }\right )^{n/d}{(1+\log _DB_0)^{n-2}\over
B_0^{d-na'(F)\over a'(F)}}D^{n+1}\sum _{l=0}^{\infty} (l+1)^{n-2}D^{-dl+nla'(F)\over a'(F)},$$
where $c_4={d\choose n}c_32^nn^{n/2}(n!)^3n^3d^{n-2}.$ 
For any $l_1\ge 0$, the lattice points $\uz\in
\Lambda$ with 
$$B_0D^{l_1+1}(m/\mofF )^{1/d}\ge \|\uz\|\ge B_0(m/\mofF )^{1/d}$$ and satisfying {\rm (1)} lie in the
union of a set of cardinality less than
$$c_5\left ({m\over\mofF }\right )^{n/d}{(1+\log _DB_0)^{n-2}\over
B_0^{d-na'(F)\over a'(F)}\det\Lambda }D^{n+1}\sum _{l=0}^{l_1}
(l+1)^{n-2}D^{-dl+nla'(F)\over a'(F)},$$
where $c_5={d\choose n}c_33^n2^{n(n-1)}n^{n/2}(n!)^3n^3d^{n-2},$ and 
less than 
$${d\choose n}n!n^3d^{n-2}(1+\log _DB_0)^{n-2}(l_1+1)^{n-1}$$ 
sublattices of smaller rank.
\endproclaim

\demo{Proof}
For a given index $l\ge 0$ let $B_l=D^lB_0$ and $C_l=DB_l$.
According to Proposition 1, the $\ux\in\br ^n$ with 
$B_l(m/\mofF )^{1/d}\le\|\ux\|\le
C_l (m/\mofF )^{1/d}$ and satisfying (1) lie in no more than
${d\choose n}\max \{n!,\ n^3\big (\log _D(B_l^{d-(n-1)a'(F)\over a'(F)}C_l^{n-1})\big )^{n-2}\}$
convex sets of the form (6) with volume no greater than
$$\multline c_32^nn^{n/2}(n!)^2 \left ({m\over\mofF }\right )^{n/d}D^n
{C_l\over B_l^{d-(n-1)a'(F)\over a'(F)}}\\
=c_32^nn^{n/2}(n!)^2
\left ({m\over \mofF }\right )^{n/d}D^{n+1}B_l^{-d+na'(F)\over a'(F)}\\
=c_32^nn^{n/2}(n!)^2
\left ({m\over \mofF }\right )^{n/d}B_0^{-d+na'(F)\over a'(F)}D^{n+1}D^{-dl+nla'(F)\over a'(F)}.
\endmultline$$

A quick estimate shows that
$$\aligned \log _D (B_l^{d-(n-1)a'(F)\over a'(F)}C_l^{n-1})&=\log _D
\big ((B_0D^l)^{d\over a'(F)}D^{n-1}\big )\\
&\le\log _D\big ((B_0D^l)^dD^{n-1}\big )\\
&< d(l+1)(1+\log _DB_0).\endaligned$$
Thus, the volume of 
all solutions $\ux$ to (1) with 
$B_0\left ({m\over\mofF}\right )^{1/d}\le\|\ux\|$
is less than
$${d\choose n}c_32^nn^{n/2}(n!)^2n^3d^{n-2}
\left ({m\over \mofF }\right )^{n/d}{(1+\log _DB_0)^{n-2}\over
B_0^{d-na'(F)\over a'(F)}}D^{n+1}\sum _{l=0}^{\infty}(l+1)^{n-2}D^{-dl-nla'(F)\over a'(F)}.$$

For the statement about the lattice points, we have fewer than
$${d\choose n}n!n^3d^{n-2}(1+\log _DB_0)^{n-2}\sum _{l=0}^{l_1}(l+1)^{n-2}\le
{d\choose n}n!n^3d^{n-2}(1+\log _DB_0)^{n-2}(l_1+1)^{n-1}$$
convex sets of the form (6), and for those containing $n$ linearly independent
lattice points, together thay  contain fewer than
$$c_5
\left ({m\over \mofF }\right )^{n/d}{(1+\log _DB_0)^{n-2}\over
B_0^{d-na'(F)\over a'(F)}\det\Lambda }D^{n+1}\sum _{l=0}^{l_1}(l+1)^{n-2}D^{-dl+nla'(F)\over a'(F)}$$
lattice points by Proposition 1.
\enddemo

\proclaim{Proposition 3} Let $F(\uX )\in\bz [\uX ]$ be a decomposable form of finite type in $n$ variables and
degree $d$. Suppose $a'(F)<d/n$. Then
for any $D\ge e$,
the integral solutions to {\rm (1)} lie in the union of a set of 
cardinality $\le c_6D^{n+1}(m/\mofF )^{n/d}$ and $\ll (1+\log _Dm+\log _D\mofF )^{n-1}$ 
proper subspaces, where $c_6=3^n2^{n(n-1)/2}n!+c_5\sum_{l=0}^{\infty} 
(l+1)^{n-2}e^{-ln/d}$. 
\endproclaim

\demo{Proof}
Let $T\in\gln (\br )$ and $S\in\gln (\bz )$ be as in the statement
of Lemma 7, and write $T=S^{-1}T'$. Consider the equivalent form
$G=F\circ S$. Then $N_G(m)=N_F(m)$, $V(G)=V(F)$ and
${\frak m}(G)=\mofF ={\Cal H}(G\circ T')$. By Lemma 7,
${\Cal H}(G)\le n^{d(n+1/2)}{\frak m}(G)^n$ and for every $\ux\in\br ^n$,
$${\frak m}(G)^{-1/d}n^{-3/2}(n!)^{-2}\| \ux\|\le\| (T')^{-1} (\ux )\|\le n^{n+1/2}{\frak m}(G)
^{(n-1)/d}\|\ux\| .$$
In other words, we may assume without loss of generality that
$$\gathered 
c_7\mofF ^{-1/d}\|\ux\|\le \| T^{-1}(\ux )\|\le c_8\mofF ^{(n-1)/d}\|\ux\|\\
\mofF ={\Cal H}(F\circ T)\\
\hofF\le c_8^d\mofF ^n,\endgathered 
\tag 8$$
where $c_7=n^{-3/2}(n!)^{-2}$ and $c_8=n^{n+1/2}.$

We will apply Proposition 2 to the lattice $\Lambda = T^{-1}(\bz ^n)$ of determinant 1,
using $B_0=1$.

Let $l_1$ be minimal such that
$$\aligned \left  (D^{l_1+1}(m/\mofF )^{1/d}\right )^{1/2}\ge\max \biggl\{ &
(c_8\mofF ^{(n-1)/d})^{1/2}
c_7^{-1}\mofF ^{1/d},\\
&\qquad (c_8\mofF ^{(n-1)/d})^{1/2}
c_8^{(d/2)}m^{1/(2d)}\mofF ^{n/2},\\
&\qquad (c_2m/\mofF )^{1/(d-na'(F))}\biggr \}.\endaligned$$
Clearly $l_1\ll 1+\log _D m+\log_D\mofF$, where the implicit constant depends only on $n$ and $d$.
Moreover, if $\| T^{-1}(\uz )\|\ge D^{l_1+1}(m/\mofF )^{1/d},$ then by
(8) we have
$$\|\uz \|^{1/2}\ge \max\{ c_7^{-1}\mofF ^{1/d}, c_8^{d/2}m^{1/(2d)}\mofF ^{n/2}\},$$
and using (8) once more,
$$m^{1/d}\hofF\le\|\uz\|\le \|T^{-1}(\uz )\|^2.\tag 9$$

By Lemma 6, (8), (9) and our choice for $l_1$, if $\uz\in\bz ^n$ is a solution to (1) with 
$\| T^{-1}(\uz )\|\ge D^{l_1+1}(m/\mofF )^{1/d},$ then there are $n$ linearly independent factors
$L_{i_1}(\uX ),\ldots ,L_{i_n}(\uX )$ of $F(\uX )$ such that
$$\aligned \left ({\prod _{j=1}^n|L_{i_j}(\uz )|\over|\det (\uL _{i_1}^{tr}\cdots 
\uL _{i_n}^{tr})|}\right )^{a'(F)}&\le {c_2m\over
\|T^{-1}(\uz ) \|^{d-na'(F)}\mofF}\\
&\le {1\over
\|T^{-1}(\uz ) \|^{(d-na'(F))/2}}\\
&\le {1\over
\|\uz  \|^{(d-na'(F))/4}}.\endaligned$$
In particular,
$${\prod _{j=1}^n|L_{i_j}(\uz )|\over|\det (\uL _{i_1}^{tr}\cdots 
\uL _{i_n}^{tr})|}< \|\uz\|^{-d+na'(F)\over 4a'(F)}<\Vert\uz\Vert ^{-1/4d}.\tag 10$$
Here we used ${d-na'(F)\over a'(F)}>1/d$.

Take such a $\uz$ and write it as a multiple of a primitive point $\uz '$;
say $\uz=g\uz ',$ for some positive integer $g$. Since $|F(\uz ')|\ge 1,$
we see that $g\le m^{1/d}$, so that $\|\uz '\|\ge\hofF$ by (9). Moreover,
we may replace $\uz$ in (10) with $\uz '$. By [E, Corollary] and
[T1, Lemma 2], such primitive $\uz'$ lie in $\ll 1$ proper subspaces.

We thus see that all lattice points $T^{-1}(\uz )$ with $\| T^{-1}(\uz )\|\ge D^{l_1+1}(m/\mofF )^{1/d}$
lie in $\ll 1$ proper subspaces. Those with $(m/\mofF )^{1/d}\le \| T^{-1}(\uz )\|\le D^{l_1+1}(m/\mofF )^{1/d}$
are dealt with by Proposition 2, using ${d-na'(F)\over a'(F)}>n/d$ again and
$D\ge e$.

It remains to deal with those lattice points with $\| T^{-1}(\uz )\|\le (m/\mofF )^{1/d}$. It is
simpler to instead estimate the number of lattice points with supnorm no greater than $(m/\mofF )^{1/d}$.
Such lattice points will be in a convex set of the form (6) with $\prod _{i=1}^na_i=(m/\mofF )^{n/d}$.
By [T1, Lemma 9], all such lattice points either lie in a proper subspace, or their number is less than
$3^n2^{n(n-1)}n!(m/\mofF )^{n/d}.$ The proposition follows.
\enddemo

Suppose $W\subseteq\br ^n$ is a subspace defined over $\bq$. Let $T_W\in\gln (\bz )$
be such that $T_W^{-1}(W)$ is the subspace spanned by the first $\dim W$ canonical basis
vectors of $\br ^n$. We will denote by $F|_W$ the decomposable form of degree $d$ in
$\dim W$ variables gotten by restricting $F\circ T_W$ to the subspace spanned by the
first $\dim W$ canonical basis vectors of $\br ^n$. Note that the solutions $\ux$
to the inequality $|F|_W(\ux )|\le m$ are in one-to-one correspondence to the
solutions to (1) lying in $W$. Further, the same holds when we consider integral
solutions.

\proclaim{Proposition 4} Suppose $F(\uX )\in\bz [\uX ]$ is a decomposable form of degree $d$ in
$n$ variables of finite type. Let $W$ be a subspace of $\br ^n$ defined over 
$\bq$ of dimension $n-1$.
There are positive constants $c_9$, $c_{10}$, $c_{11}$ and $c_{12},$
depending only on $n$ and $d$, such
that $N_{F|_W}(m)\le c_9 m^{(n-1)/d}$, and if 
${\Cal M}(F|_W)\le m^{c_{10}}$, then  
$$c_{11}m^{(n-1)/d}V(F|_W)\le N_{F|_W}(m)\le c_{12}m^{(n-1)/d}V(F|_W).$$
\endproclaim

\demo{Proof} This follows directly from
[T1, Theorem 3], Lemma 7 and [T2, Theorem 2] applied to the form $F|_W$ when
$n>2$, i.e., when $n-1\ge 2$. In the case $n=2$, the form $F|_W$ is a form
in one variable and the result is trivially valid.
\enddemo

\proclaim{Corollary}
Suppose $F(\uX )\in\bz [\uX ]$ is a decomposable form of degree $d$ in $n$
variables of finite type and $a'(F)<d/n$. 
If $\mofF \ge m^{2/n}$, then $N_F(m)\ll m^{(n-1)/d}$.
\endproclaim

\demo{Proof}
Set $D=\mofF ^{n\over 2d(n+1)}$ in Proposition 3. The
integral solutions to (1) then lie in a set of cardinality no greater than
$$c_6D^{n+1}(m/\mofF )^{n/d}=c_6(m/\mofF ^{1/2})^{n/d}\le c_6 m^{(n-1)/d}$$
and $\ll 1$ proper subspaces. By Proposition 4, the total number of solutions
in all these proper subspaces is $\ll m^{(n-1)/d}$.
\enddemo

\proclaim{Proposition 5} Suppose $F(\uX )\in\bz [\uX ]$ is a decomposable form of degree $d$ in
$n$ variables of finite type. Let $B,m\ge 1$. Suppose
 $W_1,\ldots ,W_N$ be distinct
subspaces of dimension $n-1$ defined over
$\bq$ satisfying ${\Cal M}(F|_{W_i})\le B$ for all $i$ and
$$\left ({m\over B^{n-1}}\right ) ^{1/d}\ge c_{13}(N-1),$$
where $c_{13}={2c_9\over c_{11}}\left ({n-1\over 2}\right )^{n-1}.$
Then the number of integer solutions $\ux\in\cup _{i=1}^NW_i$ to {\rm (1)} is at least
$(1/2)\sum _{i=1}^NN_{F|_{W_i}}(m).$ In particular, we have
$$2N_F(m)\ge \sum _{i=1}^N N_{F|_{W_i}}(m).$$
\endproclaim

\demo{Proof}
A simple induction argument on $N$ shows that the number of integer solutions
$\ux\in\cup _{i=1}^NW_i$ to (1) is at least
$$\sum _{i=1}^NN_{F|_{W_i}}(m)-\sum _{i=1}^{N-1}\left (
\sum _{j=i+1}^N N_{F|_{W_i\cap W_j}}(m)\right ).$$
By Proposition 4 (since $\dim (W_i\cap W_j)=n-2$),
$$\sum _{j=i+1}^N N_{F|_{W_i\cap W_j}}(m)\le
(N-1)c_{9}m^{(n-2)/d}$$
for all $i$, so that the number of integer solutions we're considering is at least
$$\sum _{i=1}^N N_{F|_{W_i}}(m)-c_9(N-1)m^{(n-2)/d}.$$

By Proposition 4 and Lemma 2 (applied to $F|_{W_i}$),
$$\aligned N_{F|_{W_i}}(m)&\ge c_{11}m^{(n-1)/d}V(F|_{W_i})\\
&\ge c_{11}m^{(n-1)/d}\left ({2\over n-1}\right )^{n-1}{\frak m}(F|_{W_i})^{-(n-1)/d}\\
&\ge c_{11}m^{(n-1)/d}\left ({2\over n-1}\right )^{n-1}{\Cal M}(F|_{W_i})^{-(n-1)/d}\\
&\ge c_{11}m^{(n-1)/d}\left ({2\over n-1}\right )^{n-1}B^{-(n-1)/d}\\
&\ge c_{11}m^{(n-2)/d}\left ({2\over n-1}\right )^{n-1}\left (m\over B^{n-1}\right )^{1/d}\\
&\ge 2c_{9}(N-1)m^{(n-2)/d}\endaligned$$
for all $i$. This together with the above estimate completes the proof.
\enddemo

\proclaim{Proposition 6} Suppose $F(\uX )\in\bz [\uX ]$ is a decomposable form of degree $d$ in
$n$ variables of finite type.
Let ${\Cal S}$ be a set of  subspaces of dimension $n-1$
defined over
$\bq$ of cardinality $N$.
Suppose $m^{1/2d}\ge\MofF ^{1/4nd}\ge c_{13}(N-1)$. Then there is a positive constant $c_{14}$,
depending only on $n$ and $d$, such that
$$\sum _{W\in{\Cal S}}N_{F|_W}(m)\ll m^{(n-1)/d}+N\big (m^{(n-1)/d}\MofF ^{-c_{14}}
+m^{(n-2)/d}(1+\log m)^{n-1}\big ).$$
\endproclaim

\demo{Proof}
Set $B=\MofF ^{c_{10}'}$, where $c_{10}'$ is the mininum of $c_{10}/2n$ and $1/4n(n-1)$.
Write
${\Cal S}$ as the union of three disjoint subsets, ${\Cal S}_1\cup {\Cal S_2}\cup{\Cal S}_3,$
where ${\Cal S}_1$ consists of the subspaces $W$ with ${\Cal M}(F|_W)\le B$,
 ${\Cal S}_2$ consists of the subspaces $W$ with
$B<{\Cal M}(F|_W)\le m^{2n}$, 
and ${\Cal S}_3$ consists of the remaining subspaces.

For the moment, set $m=\MofF ^{1/2n}$. Then since $B\le \MofF ^{1/4n(n-1)},$
$(m/B^{n-1})^{1/d}\ge \MofF ^{1/4nd}\ge c_{13}(N-1)$. 
But $B\le \MofF ^{c_{10}/2n}\le m^{c_{10}}$ also, so by Propositions 4 and 5  
$$2N_F(m)\ge \sum _{W\in{\Cal S}_1}N_{F|_W}(m)\ge c_{11} m^{(n-1)/d}\sum _{W\in{\Cal S}_1}V(F|_W)
.$$
But by [T2, Theorem 4], the integer solutions to (1) lie in the union of $\ll 1$ proper 
subspaces if $m\le\MofF ^{1/2n}$. By Proposition 4, there are no more than $c_9m^{(n-1)/d}$ solutions in a proper
subspace, so that $N_F(m)\ll m^{(n-1)/d}$. Thus,
$$\sum _{W\in{\Cal S}_1}V(F|_W)\ll 1.$$

Now let $m$ be as in the statement of Proposition 6.  
Since $B\le \MofF ^{c_{10}/2n}\le m^{c_{10}},$  
Proposition 4 and the above inequality give
$$\sum _{W\in {\Cal S}_1}N_{F|_W}(m)\le c_{12}m^{(n-1)/d}
\sum _{W\in {\Cal S}_1}V(F|_W)\ll
m^{(n-1)/d}.\tag 11$$

By [T2, Theorem 3], for any $W\in {\Cal S}_2$ we have
$$\aligned N_{F|_W}(m)&\ll 
m^{(n-1)/d}{\Cal M}(F|_W)^{-1/d}\big (1+(\log {\Cal M}(F|_W))^{n-1}
\big )\\
&\qquad +m^{(n-2)/d}\big (1+(\log m)^{n-1}+(\log {\Cal M}(F|_W))^{n-1}\big )\\
&\ll m^{(n-1)/d}B ^{-1/2d}+m^{(n-2)/d}(1+\log m)^{n-1}.\endaligned$$
Thus,
$$\sum _{W\in{\Cal S}_2}N_{F|_W}(m)\ll N
\left (m^{(n-1)/d}B ^{-1/2d}+m^{(n-2)/d}(1+\log m)^{n-1}\right ).\tag 12$$

Finally, $N_{F|_W}(m)\ll m^{(n-2)/d}$ for all $W\in {\Cal S}_3$ by [T2, Theorem 4] and
Proposition 4. Thus
$$\sum _{W\in{\Cal S}_3}N_{F|_W}(m)\ll N m^{(n-2)/d}.\tag 13$$

Proposition 6 follows from (11)-(13), setting $c_{14}=c_{10}'/2d$.

\enddemo

\head 3. Proof of Theorems 1, 2 and 3\endhead

\demo{Proof of Theorem 1}
The lower bounds for $V(F)$ and $\mofF$ in Theorem 1 are contained
in Lemmas 2 and 7.

Suppose $a'(F)<d/n$. By Lemma 1, we may assume without loss of generality that
$\hofF =\mofF =1.$ Set $D=e$ and $m=B_0=1$ in Proposition 2. Then we see that
the volume of all solutions $\ux\in\br ^n$ to (1) with
$\|\ux\|\ge 1$ is less than
$$c_4e^{n+1}\sum _{l=0}^{\infty }(l+1)^{n-2}e^{-ld+lna'(F)\over a'(F)}\ll 1.$$
Of course, the set of all $\ux\in\br ^n$ with
$\|\ux\|\le 1$ is no more than
the volume of the unit ball in $\br ^n$. 
This shows that $V(F)\ll m^{-n/d}$ when $a'(F)<d/n$, and completes the proof of
Theorem 1.
\enddemo

\demo{Proof of Theorem 2}
As in the proof of Proposition 3, we will assume (8).
Fix a $B_0\ge 1$. 
Let $S_0$ denote the cardinality of the set of integral solutions to (1) with
supnorm no greater than $c_7^{-1}m^{1/d}B_0$ and let $V_0$ denote the volume of
all real solutions to (1) with supnorm no greater than $c_7^{-1}m^{1/d}B_0$.
According to [T1, Lemma 14], 
$$|S_0-V_0|\ll 1+(m^{1/d}B_0)^{n-1}.\tag 14$$

By (8), if $\ux\in\br ^n$ is a solution to (1) with (sup)norm at least
$c_7^{-1}m^{1/d}B_0$, then $T^{-1}(\ux )$ is a solution to $|F\circ T(\uX )|\le m$
and $\| T^{-1}(\ux )\|\ge B_0(m/\mofF )^{1/d}$. Setting $D=e$ in Proposition 2, we see that the
total volume of all such $T^{-1}(\ux )$ is  less than
$$\aligned c_4\left ( {m\over \mofF}\right )^{n/d}{(1+\log B_0)^{n-2}\over 
B_0^{d-na'(F)\over a'(F)}}e^{n+1}&\sum _{l=0}^{\infty}
(l+1)^{n-2}e^{-ld+lna'(F)\over a'(F)}\\
&\ll \left ({m\over \mofF}\right )^{n/d}{(1+\log B_0)^{n-2}\over 
B_0^{d-na'(F)\over a'(F)}}.\endaligned$$
In other words,
$$m^{n/d}V(F)-V_0\ll\left ({m\over\mofF}\right )^{n/d}{(1+\log B_0)^{n-2}
\over B_0^{d-na'(F)\over a'(F)}}.\tag 15$$

Similar to the proof of Proposition 3, let $l_1$ be minimal such that
$$\aligned \left  (e^{l_1+1}(m/\mofF )^{1/d}\right )^{1/2}\ge\max \biggl\{ &
(c_8\mofF ^{(n-1)/d})^{1/2}
c_7^{-1}\mofF ^{1/d},\\
&\qquad (c_8\mofF ^{(n-1)/d})^{1/2}
c_8^{(d/2)}m^{1/(2d)}\mofF ^{n/2},\\
&\qquad (c_2m/\mofF )^{1/(d-na'(F))} \biggr \}.\endaligned$$
Then $l_1\ll 1+\log m +\log \mofF\ll 1+\log m$ (since $\mofF \le m^{1/n}$).
As in the proof of Proposition 3, if 
$$\| T^{-1}(\uz )\|\ge B_0e^{l_1+1}(m/\mofF )^{1/d}\ge
e^{l_1+1}(m/\mofF )^{1/d},$$
then we have (9) again. Arguing exactly as in the proof of Proposition
3, the integer solutions $\uz\in\bz ^n$ to (1) with $\| T^{-1}(\uz )\|\ge B_0e^{l_1+1}(m/\mofF )^{1/d}$
lie in the union of $\ll 1$ proper subspaces. By Proposition 2, then, the set of integer solutions
to (1) with (sup)norm greater than $c_7^{-1}m^{1/d}B_0$ lie in the union of a set of cardinality
$S_1$ satisfying
$$\aligned S_1< c_5\left ({m\over\mofF}\right )^{n/d}{(1+\log B_0)^{n-2}
\over B_0^{d-na'(F)\over a'(F)}}&e^{n+1}
\sum _{l=0}^{l_1}(l+1)^{n-2}
e^{-ld+lna'(F)\over a'(F)}\\
&\ll \left ({m\over\mofF}\right )^{n/d}{(1+\log B_0)^{n-2}\over 
B_0^{d-na'(F)\over a'(F)}}\endaligned\tag 16$$
and $\ll (1+\log B_0)^{n-2}(1+\log m)^{n-1}$ proper subspaces.

We now choose $B_0$. Set
$$(B_0m^{1/d})^{n-1}=\left ({m\over\mofF }\right )^{n/d}B_0^{-d+na'(F)\over a'(F)}.$$ Then
$$\gathered
B_0^{d-a'(F)\over a'(F)}=\left ({m\over \mofF ^n}\right )^{1/d}\\
B_0=\left ({m\over \mofF ^n}\right )^{a'(F)\over d(d-a'(F))}\\
B_0m^{1/d}={m^{1\over d-a'(F)}\over \mofF ^{na'(F)\over d(d-a'(F))}}\\
(B_0m^{1/d})^{n-1}={m^{n-1\over d-a'(F)}\over \mofF ^{n(n-1)a'(F)\over d(d-a'(F))}}.\endgathered\tag 17$$
Note that, by the second equation in (15) and since $\mofF\le m^{1/n}$, we indeed have $B_0\ge 1.$

Consider
the $\ll (1+\log B_0)^{n-2}(1+\log m)^{n-1}\ll (1+\log m)^{2n-3}$ proper 
subspaces above;
let ${\Cal S}$ denote this collection of proper subspaces and let $N$ denote
its cardinality. Without loss of generality, we
may assume $\dim W=n-1$ for all $W\in{\Cal S}.$  

Suppose first that $\MofF \ge m^{1/4n}$. If $m^{1/2d}\ge\MofF ^{1/4nd}\ge c_{13}(N-1),$
then by Proposition 6
$$\aligned \sum _{W\in {\Cal S}}N_{F|_W}(m)&\ll m^{(n-1)/d}+N\big (
m^{(n-1)/d}\MofF ^{-c_{14}}+m^{(n-2)/d}(1+\log m)^{n-1}\big )\\
&\ll m^{(n-1)/d},\endaligned$$
since $N\ll (1+\log m)^{2n-3}$.
If $m^{1/2d}<\MofF ^{1/4nd}$, then by Lemma 7 
$$m<\MofF ^{1/2n}\ll\mofF ^{1/2}\le m^{1/2n},$$
so $m\ll 1$. Also, if $\MofF ^{1/4nd}<c_{13}(N-1)$, then $m^{1/d}\ll (1+\log m)^{2n-3}$ and
$m\ll 1$ again. Moreover, if
$m\ll 1$, then $N\ll 1$ and by Proposition 4
$$\sum _{W\in {\Cal S}}N_{F|_W}(m)\le c_{9}N m^{(n-1)/d}\ll m^{(n-1)/d}.$$

Now suppose $\MofF < m^{1/4n}$. Then $\mofF < m^{1/4n},$ too. In this case
$$\left ({m\over \mofF ^{na'(F)/d}}\right )^{n-1\over d-a'(F)}
>m^{(n-1)/d}m^{3a'(F)(n-1)\over 4d(d-a'(F))}\gg m^{(n-1)/d}(1+\log m)^{2n-3}. $$
Again by Proposition 4,
$$\sum _{W\in {\Cal S}}N_{F|_W}(m)\le c_{9}N m^{(n-1)/d}\ll m^{(n-1)/d}(1+\log m)^{2n-3}
.$$

Thus, in all cases
$$\sum _{W\in {\Cal S}}N_{F|_W}(m)\ll
\left ({m\over \mofF ^{na'(F)/d}}\right )^{n-1\over d-a'(F)}
(1+\log m)^{n-2}.\tag 18$$
Theorem 2 follows from (14)-(18).

\enddemo

\demo{Proof of Theorem 3}
By the corollary to Propositions 3 and 4, we only need to deal with
the case where
$\mofF\le m^{2/n}$. Set $D=e$ in Proposition 3. 
Since 
$\log\mofF\le (2/n)\log m$,
Proposition 3 shows that the integral solutions to (1) lie in the union of
a set of cardinality $\ll (m/\mofF )^{n/d}$ and $\ll (1+\log m)^{n-1}$ proper
subspaces which we may assume, without loss of generality, are all dimension
$n-1$. 

Now if $\MofF \ge m^{1/4n}$, then we argue exactly as in the proof of 
Theorem 2 above and conclude that the number
of solutions contained in these proper subspaces is $\ll m^{(n-1)/d}.$

Suppose $\MofF < m^{1/4n}$. By Proposition 4, the number of solutions
contained in these proper subspaces is $\ll m^{(n-1)/d}(1+\log m)^{n-1}$.
But $$m^{(n-1)/d}(1+\log m)^{n-1}\ll {m^{n/d}\over m^{1/4d}}\le \left ({
m\over\MofF}\right )^{n/d}\le\left ({m\over\mofF }\right )^{n/d}.$$

\enddemo

\head 4. Proof of Theorems 4 and 5\endhead

We need some notation from [T2]. Let $F(\uX )=\prod _{i=1}^dL_i(\uX )$ be a factorization
of $F$ as in the statement of Lemma 3. Consider the $d^n$-dimensional vector
with components $\det (\uL _{i_1}^{tr}\cdots\uL _{i_n}^{tr})$. The quantity $Q(F)$ is defined to be
the infimum over all such factorizations of the $L^2$ norms of these vectors.
Let
$$NS(F)=\prod _{(i_1,\ldots ,i_n)}'{|\det (\uL _{i_1}^{tr}\cdots \uL
_{i_n}^{tr})|\over \Vert \uL _{i_1}\Vert\cdots\Vert \uL _{i_n}\Vert}.$$
Here the restricted product is over those $(i_1,\ldots ,i_n)$ where
$\uL _{i_1},\ldots ,\uL _{i_n}$ are linearly independent.

\demo{Proof of Theorem 4}
Lemma 3 shows that 
$$Q(F)^2\ge (d/n)^nn!\mofF ^{2n/d}.\tag 19$$
If $F(\uX )\in\bz [\uX ]$, then $|NS(F)|^{-1}\le\hofF ^{d\choose n}$
by [T1, Lemma 3].
If we assume $\hofF =\MofF$ and further that $F$ doesn't vanish at any non-trivial rational
point, then this together with Lemma 7 shows that
$$|\log NS(F)|\ll |\log\big (\mofF\big )|.\tag 20$$
By [T2, Thereom 1],
$$V(F)\ll Q(F)^{-1}(1+|\log NS(F)|)^{n-1}.$$
Theorem 4 follows from this, (19) and (20).
\enddemo

To prove Theorem 5, we note that if $a'(F)=d/n$ in the proof of Proposition 1,
then (7) becomes
$${\prod _{j=1}^n|L_{i_j}(\ux )|\over |\det (\uL _{i_1}^{tr}\cdots\uL _{i_n}^{tr})|}\le c_3,
\tag 7'$$
Moreover, the hypothesis $\hofF =\mofF$ used to obtain (7) is not 
necessary here. Thus, the hypothesis $\hofF =\mofF$ in Proposition 1 is
unnecessary when $a'(F)=d/n$.

\demo{Proof of Theorem 5}
By [T2, Lemma 4], any solution $\ux$ to (1) satisfies an inequality of the form
$${\prod _{j=1}^n|L_{i_j}(\ux )|\over|\det (\uL _{i_1}^{tr}\cdots \uL _{i_n}^{tr})|}\le
\left ({c_{15}m\over \|\ux\| ^{d-na(F)}\hofF NS(F)^{d-(n-1)a(F)}}\right )^{1/a(F)},\tag 21$$
where $c_{15}$ is a constant which depends only on $n$ and $d$ and $a(F)$ is the same as in [T1, Theorem 3] above.

Choose $l_1$ minimal such that
$$\big (e^{l_1+1}(m/\mofF )^{1/d}\big )^{(d-na(F))/2} \ge\max \left \{ \big (m^{1/d}\hofF\big )^{(d-na(F))/2},
{c_{15}m\over \hofF NS(F)^{d-(n-1)a(F)}}\right \}.$$
By (8) and (20), we see that $l_1\ll 1+\log\mofF$ if 
$\mofF ^{-n/d}(1+\log\mofF )^{n-1}\le m^{-1/d}$. 
For any solution  $\ux\in\br ^n$ to (1) which satisfies $\|\ux\|> e^{l_1+1}(m/\mofF )^{1/d},$
we have
$${\prod _{j=1}^n|L_{i_j}(\ux )|\over|\det (\uL _{i_1}^{tr}\cdots \uL _{i_n}^{tr})|}<
\|\ux\| ^{-(d-na(F))\over 2a(F)}$$
and also $\|\ux\|\ge m^{1/d}\hofF$. As in the proof of Proposition 3, such $\ux\in\bz ^n$ lie in
$\ll 1$ proper subspaces. 

For an index $l\ge 0,$ set $D=e$, $B=e^l$ and $C=eB$ in Proposition 1. Then
$${d\choose n}\max\{ n!,n^3\big (\log _D(BC^{n-1})\big )^{n-2}\}<{d\choose n}n^{n+1}(1+l)^{n-2}.$$
All integral solutions $\ux$ to (1) with $(m/\mofF )^{1/d}e^l\le\|\ux\|\le (m/\mofF )^{1/d}e^{l+1}$
thus lie in the union of a set of cardinality $\ll (m/\mofF )^{n/d}(1+l)^{n-2}$ and $\ll (1+l)^{n-2}$
proper subspaces. 

There are $\ll (m/\mofF )^{n/d}$  $\ux\in\bz ^n$ with $\|\ux\|\le (m/\mofF )^{1/d}.$ Since
$\mofF \ge 1$ we have $\mofF ^{-n/d}\le m^{-n/d}(1+\log\mofF )^{n-1}\le m^{-1/d}$, so that there are
$\ll m^{(n-1)/d}$ solutions $\ux\in\bz ^n$ to (1) with $\|\ux\|\le (m/\mofF )^{1/d}$. Also,
$$\aligned \left ({m\over\mofF }\right )^{n/d}\sum _{l=0}^{l_1}(1+l)^{n-2}&\le 
\left ({m\over\mofF }\right )^{n/d}(1+l_1)^{n-1}\\
&\ll 
\left ({m\over\mofF }\right )^{n/d}(1+\log\mofF )^{n-1}\\
&\le m^{(n-1)/d}.\endaligned$$
In this way we see that the solutions $\ux\in\bz ^n$ to (1) with $\|\ux\|\le e^{l_1+1}(m/\mofF )^{1/d}$ 
lie in the union of a set of cardinality
$\ll m^{(n-1)/d}$ and $\ll (1+\log\mofF )^{n-1}$ proper subspaces. The larger solutions lie in
$\ll 1$ proper subspaces. Denote the total number of proper subspaces here by $N$. As before,
we may assume all are dimension $n-1$.

Suppose $m\ge \MofF ^{1/2}$. If $\MofF ^{1/4nd}\ge c_{13}(N-1)$, then by Proposition 6 the
total number of solutions in our subspaces is $\ll m^{(n-1)/d}$. If $\MofF ^{1/4nd}<
c_{13}(N-1)$, then $\MofF \ll 1$ and we have $\ll 1$ proper subspaces. In this case,
Proposition 4 implies that the total number of solutions in these subspaces is
$\ll m^{(n-1)/d}$.

It remains to deal with the case where $\MofF > m^2$. In this case, Lemma 7
shows that $\mofF\gg m^{2/n}$; we set $D=\mofF ^{n\over 2d(n+1)}$
and choose $l_1$ minimal such that
$$\big (D^{l_1+1}(m/\mofF )^{1/d}\big )^{(d-na(F))/2} \ge\max \left \{ \big (m^{1/d}\hofF\big )^{(d-na(F))/2},
{c_{15}m\over \hofF NS(F)^{d-(n-1)a(F)}}\right \}.$$
Now $l_1\ll 1$ and, as above, the solutions $\ux\in\bz ^n$ with $\|\ux\|\ge D^{l_1+1}(m/\mofF )^{n/d}$ lie
in $\ll 1$ proper subspaces. 

For an index $l\ge 0$ let $B_l=D^l$ and $C_l=D^{l+1}$ in Proposition 1. The solutions $\ux\in\bz ^n$ to (1) with
$(m/\mofF )^{1/d}B_l\le\|\ux\|\le (m/\mofF )^{1/d}C_l$ lie in the union of $\ll (1+l)^{n-1}$ proper subspaces
and a set of cardinality no greater than
$$\aligned c_33^n2^{n(n-1)/2}n^{n/2}(n!)^2\left (m\over\mofF \right )^{n/d}D^n{C_l\over B_l}&\ll
\left (m\over\mofF \right )^{n/d}D^{n+1}\\
&={m^{n/d}\over\mofF ^{n/2d}}\\
&\ll m^{(n-1)/d}.\endaligned$$
Taking the sum from $l=0$ up to $l=l_1\ll 1$, we see that all integral solutions $\ux$ to (1)
with $(m/\mofF )^{1/d}D^{l_1+1}\ge \|\ux\|\ge (m/\mofF )^{1/d}$ lie in 
the union of a set of cardinality $\ll m^{(n-1)/d}$ and $\ll 1$ proper subspaces. Larger solutions
lie in the union of $\ll 1$ proper subspaces, and there are $\ll 1+(m/\mofF )^{n/d}\ll m^{(n-1)/d}$
$\ux\in\bz ^n$ with $\|\ux\|\le (m/\mofF )^{1/d}$. 
Proposition 4 applied to the subspaces completes the proof of this case, and 
thus the proof of Theorem 5.
\enddemo

\head 5. Some Examples\endhead

Fix an even $d\ge 4$ and $0<\epsilon \le 1/3$. Let 
$$F_{\epsilon}(X,Y)=\big (X^l-(\epsilon Y)^l\big )\big ((\epsilon X)^l-Y^l\big ),$$
where $l=d/2$. Then
$$F_{\epsilon}(X,Y)=\prod _{i=1}^l(X-\rho ^i\epsilon Y)(\rho ^i\epsilon X-Y),$$
where $\rho$ is a primitive $l$th root of unity. 

Suppose $1\le y\le (3\epsilon )^{-1/2}$ and $\epsilon y\le x\le 1/(3y)$. Then 
$$\gathered |x-\rho ^i\epsilon y|\le x+\epsilon y\le 2x\\
|\rho ^i\epsilon x-y|\le \epsilon x+y\le 1/(9y)+y\le 10y/9 \endgathered$$
for any $i$.
In particular, $|F_{\epsilon }(x,y)|\le (2x)^l(10y/9)^l\le (2/3)^l(10/9)^l<1.$ From this, we
see that
$$\aligned V(F_{\epsilon })&>\int _1^{(3\epsilon )^{-1/2}}\int _{\epsilon y}^{(3y)^{-1}}dx\, dy\\
&=\int _1^{(3\epsilon )^{-1/2}}(3y)^{-1}-\epsilon y\, dy\\
&>{-\log 3-\log \epsilon -1\over 6}.\endaligned$$
Thus, if $\epsilon <(3e)^{-2}$ we have 
$$V(F_{\epsilon })>{-\log\epsilon \over 12}.\tag 22$$

Now suppose $T\in\gln (\br )$ with $|\det T|=1$ and write $F_{\epsilon }(X,Y)=\prod _{i=1}^dL_i(X,Y)$,
where 
$$L_i(X,Y)=\cases X-\rho ^i\epsilon Y&\text{if $i\le l$,}\\
\rho ^{d-i}\epsilon X -Y&\text{if $l<i\le d$.}\endcases$$
By Hadamard's inequality, for $1\le i\le l$ we have
$$\|\uL_i T\|\cdot\|\uL _{i+l}T\|\ge |\det (\uL _i^{tr}\uL _{i+l}^{tr})|=1-\epsilon ^2$$
since $|\det T|=1$. Thus, ${\Cal H}(F_{\epsilon }\circ T)\ge (1-\epsilon ^2)^l$ and
$$(1+\epsilon ^2)^l={\Cal H}(F_{\epsilon })\ge {\frak m}(F_{\epsilon })\ge (1-\epsilon ^2)^l.\tag 23$$

Since we may choose $\epsilon$ arbitrarily small, (22) and (23) show that $V(F)$ cannot be
bounded above by a function of $\mofF$ in the case where $d$ is even and $n=2$.

Now let $\epsilon =p^{-1/l}$ for a large prime $p$, for example. Then $p^2F_{\epsilon }(X,Y)
=(pX^l-Y^l)(X^l-pY^l)\in\bz [X,Y]$
is finite type, even. Moreover, by Lemma 1, (22) and (23), we have
$$\gathered p^{-4/d}\log p\ll V(p^2F_{\epsilon })\\
p^{-4/d}\gg\ll {\frak m}(p^2F_{\epsilon })^{-2/d},\endgathered$$
with absolute implicit constants.
In particular, 
$$V(p^2F_{\epsilon })\gg {\frak m}(p^2F_{\epsilon })^{-2/d}\log {\frak m}(p^2F_{\epsilon }).$$
This shows that Theorem 4 is best possible.

\enddocument